\documentclass{elsart3-1}

\usepackage{graphicx}
\usepackage{epstopdf}

\usepackage{amssymb}
\usepackage{amsmath}

\usepackage[english,francais]{babel}

\newtheorem{theorem}{Theorem}[section]
\newtheorem{conjecture}{Conjecture}[section]
\newtheorem{lemma}[theorem]{Lemma}
\newtheorem{e-proposition}[theorem]{Proposition}

\newtheorem{e-definition}[theorem]{Definition\rm}

\renewcommand{\theequation}{\arabic{equation}}
\def\og{\leavevmode\raise.3ex\hbox{$\scriptscriptstyle\langle\!\langle$~}}
\def\fg{\leavevmode\raise.3ex\hbox{~$\!\scriptscriptstyle\,\rangle\!\rangle$}}

\begin{document}
Number theory
\centerline{}
\begin{frontmatter}


\thanks[label1]{Partially supported by Serbian Ministry of Education and Science, Project 174032.
}
\selectlanguage{english}
\title{The number of unimodular roots of some reciprocal
polynomials}


\selectlanguage{english}
\author{Dragan Stankov},
\ead{dstankov@rgf.bg.ac.rs}

\address{
Katedra Matematike RGF-a,
University of Belgrade,
Belgrade, \DJ u\v sina 7,
Serbia}



\begin{abstract}
\selectlanguage{english}

We introduce a sequence $P_{2n}$ of monic reciprocal polynomials with integer coefficients having the central coefficients fixed. We prove that the ratio between number of
nonunimodular roots of $P_{2n}$ and its degree $d$ has a limit when $d$ tends to infinity. We present an algorithm for calculation the limit and a numerical method for its approximation. If $P_{2n}$ is the sum of a fixed number of monomials we determine the central coefficients such that the ratio has the minimal limit. We generalise the limit of the ratio for multivariate polynomials. Some examples suggest a theorem for polynomials in two variables which is analogous to Boyd's limit formula for Mahler measure.

\vskip 0.5\baselineskip



\end{abstract}
\end{frontmatter}


\selectlanguage{english}
\section{Introduction}
If $P(x)=a_dx^d+a_{d-1}x^{d-1}+\cdots +a_1x+a_0$ ($a_d \ne 0$) has zeros $\alpha_1,\alpha_2,\ldots,\alpha_d$ then the Mahler measure of $P(x)$ is \[M(P(x)) = |a_d|\prod_{j=1}^{d}
\max(1, |\alpha_j|).\]
Let $I(P)$ denote the the number of complex
zeros of $P(x)$ which are $<1$ in modulus, counted with multiplicities.
Let $U(P)$ denote the number of zeros of $P(x)$ which are $=1$ in modulus,
(again, counting with multiplicities). Such zeros 
are called unimodular.
Let $E(P)$ denote the number of complex
zeros of $P(x)$ which are $>1$ in modulus,, counted with multiplicities.
Then it is obviously that $I(P)+U(P)+E(P)=d$.
Pisot number can be defined as a real algebraic integer greater than 1 having the minimal polynomial $P(x)$ of degree $d$ such that $I(P)=d-1$.
Salem number 
is a real algebraic integer $> 1$ having the minimal polynomial $P(x)$ of degree $d$ such that $U(P)=d-2$, $I(P)=1$.

We say that a polynomial of degree $d$ is reciprocal if
$P(x) = x^dP(1/x)$.
If moduli of coefficients are small then a reciprocal polynomial has many unimodular roots. A Littlewood polynomial is a polynomial all of whose coefficients are $1$ or $-1$. Mukunda \cite{Muk}
showed that every self-reciprocal Littlewood polynomial of odd
degree at least 3 has at least 3 zeros on the unit circle. Drungilas \cite{Dru} proved that every
self-reciprocal Littlewood polynomial of odd degree $n \ge 7$ has at least $5$ zeros on the unit
circle and every self-reciprocal Littlewood polynomial of even degree $n \ge 14$ has at least
$4$ unimodular zeros. In \cite{BCFJ} two types of very special Littlewood polynomials are considered: Littlewood polynomials with one sign change in the sequence of coefficients and Littlewood
polynomials with one negative coefficient. The numbers $U(P)$ and $I(P)$ of such Littlewood
polynomials $P$ are investigated. In \cite{BEFL} Borwein, Erd\'{e}lyi, Ferguson and Lockhart showed that there exists a cosine polynomials $\sum_{m=1}^{N}\cos(n_m\theta) $
with the $n_m$ integral and all different so that the number of its real zeros in $[0,2\pi)$ is
$O(N^{9/10}(\log N)^{1/5})$ (here the frequencies $n_m = n_m(N)$ may vary with $N$). However, there
are reasons to believe that a cosine polynomial $\sum_{m=1}^{N}\cos(n_m\theta) $ always has many zeros in
the period.

Clearly,
if $\alpha_j$, is a root of a reciprocal $P(x)$ then $1/\alpha_j$ is also a root of $P(x)$ so that $I(P)=E(P)$.
Let $n$, $k$, $a_0,a_1,\ldots,a_k$, be integers such that $n>k\ge 0$, and let $P_{2n}(x)$ be a monic, reciprocal polynomial with integer coefficients
$$P_{2n}=x^n\left(x^n+a_0+\frac{1}{x^n}+\sum
_{j=1}^k a_j\left(x^j+\frac{1}{x^j}\right)\right).$$
Let $C(P)=\frac{I(P)+E(P)}{2n}$ be the ratio between number of nonunimodular zeros of $P$ and its degree. Actually, it is the probability that a randomly chosen zero is not unimodular, and $C(P)=\frac{E(P)}{n}$.

\section{The main theorem}

\begin{theorem}
If $k>0$ is an integer then for all fixed integers $a_j$, $j=1, \ldots,k$ there is a limit $C(P_{2n})$ when $n$ tends to infinity.
\end{theorem}
\begin{pf}
The theorem will be proved if we show that $1-C(P_{2n})$  has a limit when $n$ tends to $\infty$. Since $1-C(P_{2n})=\frac{U(P_{2n})}{2n}$ we have to count the unimodular roots of $P_{2n}(x)$.
If we use the substitution $x=e^{it}$ in the equation $P_{2n}(x)=0$ we get
$$e^{int}\left(2\cos nt +a_0+\sum_{j=1}^k 2 a_j\cos jt\right)=0$$
Since $e^{int}\ne 0$ it follows that the equation is equivalent to
\begin{equation}\label{eq:CES}
\cos nt =-\frac{a_0}{2}-\sum_{j=1}^k a_j\cos jt.
\end{equation}
From the substitution $x=e^{it}$ it follows that $x$ is unimodular if and only if $t$ is real so that we have to count the real roots of \eqref{eq:CES} ($t\in [0,2\pi))$.
If $\Gamma_1$ is the graph of $f_1(t)=\cos nt$ 
and $\Gamma_2$ is the graph of $f_2(t)=-a_0/2-\sum_{j=1}^k a_j\cos jt$, the function on the right side of equation \eqref{eq:CES}, then $U(P)$ is equal to the number of intersection points of these two graphs.
These intersection points are obviously settled between lines $y=-1$ and $y=1$.
Graph $\Gamma_2$ of the continuous function $f_2$ is fixed i.e. does not depend on $n$, therefore we can introduce a partition of $[0,2\pi]$ using points $0=t_0<t_1<\ldots<t_p=2\pi$ such that $|f_2(t_j)|=1$, $0<j<p$. Let us consider subintervals $I_j=[t_{j-1},t_{j}]$ such that if $t\in I_j$ then $|f_2(t)|<1$, $j\in J=\{j_1,j_2,\ldots,j_{r}\}\subseteq \{1,2,\ldots,p\}$.
\begin{e-definition}
A part of the graph of $f_1(t)=\cos nt$ such that $(k-1)\pi/n \le t \le k\pi/n$, $k\in \mathbb{Z}$ is $k$-th branch of $\cos nt$. The interval $[(k-1)\pi/n,k\pi/n]$ is the domain of the $k$-th branch.
\end{e-definition}
Each branch of $\cos nt$ obviously has exactly one intersection point with the $t$-axis. We are going to prove that if $n$ is large enough then each branch of $\cos nt$ also has
exactly one intersection point with $\Gamma_2$.
We need the next lemma which will be proved in the next subsection.

\begin{lemma}\label{lemmaBeps}
For all $B_1, B_2>0$ and $ \varepsilon$ such that $1> \varepsilon>0$, there is $n_0\in \mathbb{N}$ such that if $n\ge n_0$ then

(1) if $|\cos(nt)|<1-\varepsilon$ then $n|\sin(nt)|>B_1$,

(2) if $|\cos(nt)|>1-\varepsilon$ then $n^2|\cos(nt)|>B_2$.
\end{lemma}

We will also need the following claims.
\begin{enumerate}
  \item There is a bound $B_1$ of the modulus of the first derivative of $f_2(t)$. Indeed
  
  \noindent $|f'_2(t)|=$ $|\sum_{j=1}^k ja_j\sin jt|\le \sum_{j=1}^k j|a_j|=:B_1$.
  \item There is a bound $B_2$ of the modulus of the second derivative of $f_2(t)$. Indeed
  
  \noindent $|f''_2(t)|=$ $|\sum_{j=1}^k j^2a_j\cos jt|\le \sum_{j=1}^k j^2|a_j|=:B_2$.

  \item The first derivative of $f_2(t)$ has a finite number of roots on $[0,2\pi]$ so that
there is $ \varepsilon_j>0$ such that $1-\varepsilon_j$ is greater than the value at each local maximum and $-1+\varepsilon_j$ is less than the value at each local minimum of $f_2(t)$ on $(t_{j-1},t_{j})$.
  \item If the domain of a branch of $\cos nt$ is the subset of the interior of $I_j$ then $\cos nt-f_2(t)$ has values of the opposite sign at the end points of the domain so that the branch has at least one intersection point with $\Gamma_2$.
\end{enumerate}

Since $f_2(t)$ is continuous at $t_{j-1}$ and $t_j$ it follows that
there are $\delta_{1j}>0$, $\delta_{2j}>0$ such that if $t\in(t_{j-1},t_{j-1}+\delta_{1j})$ or $t\in(t_{j}-\delta_{2j},t_{j})$ then $1-|f_2(t)|<\varepsilon_j$. If we bring to mind (iii) it follows that $f_2(t)$ is monotonic on $(t_{j-1},t_{j-1}+\delta_{1j})$ and on $(t_{j}-\delta_{2j},t_{j})$. Therefore we can choose $\delta_{1j}>0$, $\delta_{2j}>0$ such that $|f_2(t_{j-1}+\delta_{1j})|=1-\varepsilon_j$, $|f_2(t_{j}-\delta_{2j})|=1-\varepsilon_j$


Using Lemma \ref{lemmaBeps} (1) there is $n_j$ such that if $n\ge n_j$ and $|f_1(t)|<1-\varepsilon_{j}$ then $|f_1'(t)|>B_1$. It follows that on $E_j:=[t_{j-1}+\delta_{1j},t_j-\delta_{2j}]$ a branch of $\cos nt$ and $\Gamma_2$ can not have more than one intersection point: if they have two intersection points $M_1$, $M_2$ then using the mean value theorem for the continuous function $f_1$ the slope $S$ of the line $M_1M_2$ is greater than $B_1$ in modulus. Using the mean value theorem again for the continuous function $f_2$ it follows that there is a point $t$ such that $f'_2(t)=S$ so that $|f'_2(t)|=|S|>B_1$ which is the contradiction with (i).

It remains to be proved that if the domain of a branch is the subset of $D_{1j}=(t_{j-1},t_{j-1}+\delta_{1j}]$ or of $D_{2j}=[t_{j}-\delta_{2j},t_{j})$ then the branch of $\cos nt$ and $\Gamma_2$ can not have more than one intersection point. Let $1>f_2(t)>1-\varepsilon_j$ and let the branch has an adjacent branch such that the union of its domains is $[(k-1)\pi/n,(k+1)\pi/n]\subset D_{1j}$ and $k$ is even. Then using Lemma \ref{lemmaBeps} (2) it follows that if $\cos(nt)>1-\varepsilon$ then $f''_1(t)-f''_2(t)=-n^2\cos nt-f''_2(t)<-B_2-f''_2(t)$ is negative. Therefore $f_1(t)-f_2(t)$ is a concave function so that its graph can have at most two intersection points with the line $y=0$. If such an adjacent branch does not exist which means that $t_{j-1}\in [k\pi/n,(k+1)\pi/n]$, $k$ is even, then we can prove the concavity of $f_1(t)-f_2(t)$ in the same manner. We conclude that if $t_{j-1}$, the start point of $I_j$, is in the domain of a branch of $\cos nt$ then the branch can have 0, 1, or 2 intersection points with $\Gamma_2$ (see fig. 1).

If $-1<f_2(t)<-1+\varepsilon_j$ after showing the convexity of $f_1(t)-f_2(t)$ on $D_1$ the claim follows in the similar manner.
Analogously we prove the claim if the domain of a branch is the subset of $D_2$ as well as the claim for the end point of $I_j$: if $t_{j}$ is in the domain of a branch of $\cos nt$ then the branch can have 0, 1, or 2 intersection points with $\Gamma_2$.
\begin{figure}[t]
  \centering
  \includegraphics[width=0.5\textwidth, height=0.5\textheight]{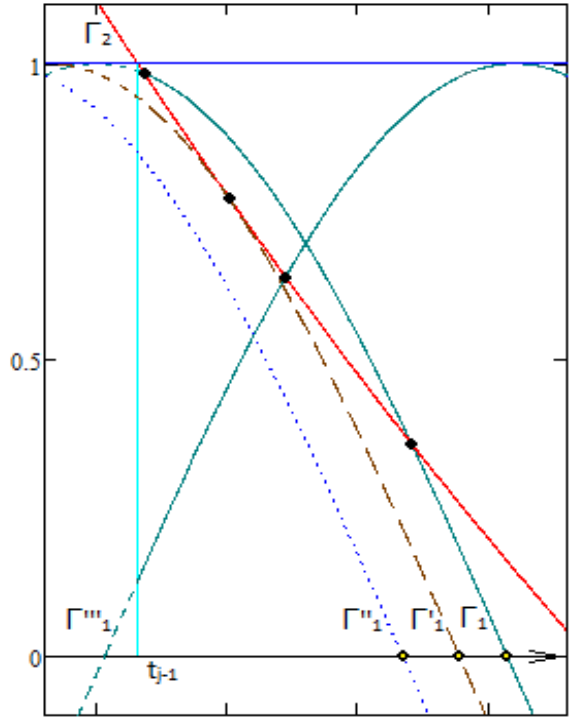}
\caption{
If $t_{j-1}$, the start point of $I_j$, is in the domain of a branch of $\cos nt$ then the branch can have 0, ($\Gamma''_1$) 1 ($\Gamma'_1$, $\Gamma'''_1$), or 2 ($\Gamma_1$) intersection points with $\Gamma_2$.}\label{fig:5}       
\end{figure}

We conclude that if $n$ is large enough then each branch of $\cos nt$, such that the start and the end point of $I_j$ are not elements of its domain, has exactly one intersection point with $\Gamma_2$.  
Thus the number $U_j$ of intersection points of $\Gamma_1$ and $\Gamma_2$ differs to the number $V_j$ of intersection points of $\Gamma_1$ and the $t$-axis, $t\in I_j$, by 0,1 or 2, because in the beginning and at the end of $I_j$ branches are not complete (see fig. 1).
If we take the sum $U_j$ and $V_j$ over all $r$ subintervals then it is clear that 
$U(P_{2n})$ differs to the number $V(P_{2n})=\sum_{j\in J} V_j$ 
by a number $\le 2r$. Since $2r$ does not depend on $n$ it follows that
$$\left(\lim_{n\rightarrow \infty}(1-C(P_{2n}))=\right)\lim_{n\rightarrow \infty}\frac{U(P_{2n})}{2n}=\lim_{n\rightarrow \infty}\frac{V(P_{2n})}{2n}\left(=\lim_{n\rightarrow \infty}\frac{\sum_{j\in J} V_j}{2n}\right).$$
Since the intersection points of the graphs of $y=\cos nt$ and the $t$-axis are obviously uniformly distributed on $I_j$ we conclude
$$\lim_{n\rightarrow \infty}\frac{\sum_{j\in J} V_j}{2n}=\frac{\sum_{j\in J} |I_j|}{2\pi}.$$

\end{pf}
\qed

\subsection{Proof of Lemma \ref{lemmaBeps}}

Using the symmetry and the periodicity of $\cos nt$ it is enough to prove the claim for the first branch of $\cos nt$, $t\in [0,\pi/n]$.
For an arbitrarily chosen $\varepsilon>0$ and $n\in N$ we determine $\tau$ such that $|\cos(n\tau)|=1-\varepsilon$.
It follows that $\tau=\arccos(1-\varepsilon)/n$ or $\tau=\arccos(-1+\varepsilon)/n$ so that

(1) if $t\in (\tau,\pi/n-\tau)$ then $n\sin nt>n\sin n\tau=n \sin(\arccos(1-\varepsilon)) \rightarrow \infty$ when
$n \rightarrow \infty$. Therefore the claim follows immediately if we chose $$n_1=\left\lceil \frac{B_1}{\sin(\arccos(1-\varepsilon))}\right\rceil .$$

(2) if $t\in (0,\tau)\bigcup(\pi/n-\tau,\pi/n)$ then $n^2|\cos nt|>n^2|\cos n\tau|=n^2 \cos(\arccos(1-\varepsilon))=n^2 (1-\varepsilon) \rightarrow \infty$ when
$n \rightarrow \infty$. Therefore the claim follows immediately if we chose $$n_2=\left\lceil \sqrt{\frac{B_2}{1-\varepsilon}}\right\rceil .$$

It remains to take $n_0=\max(n_1,n_2)$.
\qed

\subsection{Algorithm for determination $\lim_{n\rightarrow\infty}C(P_{2n})$}

In the proof of Theorem 1 we actually declared steps of an algorithm for determination $\lim_{n\rightarrow\infty}C(P_{2n})$:
\begin{enumerate}
\item determine all real roots $t_j$ of the equations $f_2(t)=1$ and $f_2(t)=-1$ ,

\item arrange them as an increasing sequence $0=t_0<t_1<\ldots<t_p=2\pi$,

\item determine $I_j=[t_{j-1},t_{j}]$ such that if $t_{j-1}<t<t_{j}$ then $|f_2(t)|<1$, $j\in J=\{j_1,j_2,\ldots,j_r\}\subseteq \{1,2,\ldots,p\}$,

\item calculate $\lim_{n\rightarrow\infty}C(P_{2n})=1-\sum_{j\in J} (t_j-t_{j-1})/(2\pi).$
\end{enumerate}

If $\overline{f_2(t)}$ is defined:

\begin{equation*}
    \overline{f_2(t)} = \begin{cases}
               1, & |f_2(t)|\ge 1\\
               0, & \text{otherwise}
           \end{cases}
\end{equation*}
then
\begin{equation}\label{eq:Cint}
\lim_{n\rightarrow\infty}C(P_{2n})=\frac{1}{2\pi}\int_{0}^{2\pi}\overline{f_2(t)}dt.
\end{equation}

\section{Approximating $\lim_{n\rightarrow\infty}C(P_{2n})$}

The equation $f_2(t)=\pm 1$ i.e. $-a_0/2-\sum_{j=1}^k a_j\cos jt=\pm 1$ is algebraic in $\cos t$ so that $t_j$ can be expressed by arccosine of an algebraic real number $\alpha \in [-1,1]$ thus only solutions of this kind should be taken into account.

We can approximate numerically the integral in \eqref{eq:Cint} i.e. $\lim_{n\rightarrow\infty}C(P_{2n})$. Suppose the interval $[0,2\pi]$ is divided into $p$ equal subintervals of length $\Delta t=2\pi/p$ so that we introduce a partition of $[0,2\pi]$ $0=t_0<t_1<\ldots<t_p=2\pi$ such that $t_{j}-t_{j-1}=\Delta t$. Then we chose numbers $\xi_j\in [t_{j},t_{j-1}]$   and count all $\xi_j$ such that $|f_2(\xi_j)|>1$, $j=1,2,\ldots,p$. If there are $s$ such $\xi_j$ then $\lim_{n\rightarrow\infty}C(P_{2n})$ is approximately equal to $\frac{s}{p}$.

$$\lim_{n\rightarrow\infty}C(P_{2n})\approx\frac{1}{p}\sum_{j=1}^p \overline{f_2(j\frac{2\pi}{p})}$$
where we chosed $\xi_j=2j\pi/p$.

\subsection{Small limit points of $C(P_{2n})$}


In the case of trinomials i.e. if  $k=0$, $|a_0|\le2$ then all roots of $P_{2n}(x)=x^{2n}+a_0x^n+1$ obviously are unimodular. If $|a_0|>2$ then $P_{2n}$ does not have any unimodular root so that $C(P)$ tends either to zero or to one as $n$ approaches infinity.

In the case of quadrinomials i.e.
if  $k=1$, $a_0=0$, $a_1=\pm 1$ then $P_{2n}(x)=x^{2n}\pm x^{n+k}\pm x^{n-k}+1=(x^{n-k}\pm 1)( x^{n+k}\pm 1)$ so that obviously all roots are unimodular.
If $|a_1|>1$ then $$C(x^{2n}+a_1 x^{n+k}+a_1 x^{n-k}+1)=2\arccos(1/a_1)/\pi$$ so that it has the minimum value 2/3 when $a_1=2$ and $C(P)$ tends to one as $a_1$ approaches infinity.

If we exclude trinomials and quadrinomials then it is clear that the limit points of $C(P_{2n})$ are always greater than zero.
A natural question that arises here is what is the smallest value, greater than 0, of the limit points of $C(P_{2n})$?

\subsection{Hexanomials with smallest limit points of $C(P_{2n})$}

Between all pentanomials $x^{2n}+a_kx^{n+k}+a_0x^n +a_kx^{n-k}+1$ an exhaustive search such that $k=1,2,\ldots,10$, $a_k=\pm1,\pm 2,\ldots \pm 10$ suggests that $C(x^{2n}+x^{n+1}+x^n +x^{n-1}+1)$ has the minimal limit point. It is equal to $$1/\pi\arccos(1/2)=1/3.$$

We have submitted an exhaustive search between all hexanomials $x^{2n}+a_{j_2}x^{n+j_2}+a_{j_1}x^{n+j_1}+a_{j_1}x^{n-j_1} +a_{j_2}x^{n-j_2}+1$ such that $j_1=1,2,\ldots,10$,  $j_2=j_1+1,j_1+2,\ldots,10$, $a_{j_1}=\pm1,\pm 2,\ldots \pm 10$,$a_{j_2}=\pm1,\pm 2,\ldots \pm 10$ suggests that $C(x^{2n}+x^{n+3}+x^{n+1} +x^{n-1}+x^{n-3}+1)$ has the minimal limit point.

Using the algorithm we solve the equation
$\cos3t+\cos t=\pm1$. Since $\cos 3t=4\cos^3 t-3\cos t$ if we substitute $\cos t=x$ we get algebraic equations $4x^3-2x=\pm1$ each of them with unique real solution $\pm \alpha$ where
$$\alpha=\frac{\sqrt[3]{\frac{\sqrt{57}}{288}+\frac{23}{864}}+\frac{1}{6}}{\sqrt[3]{\frac{\sqrt{57}}{72}+\frac{1}{8}}}\approx 0.885,\;\;\;C(P_{2n})\rightarrow \frac{2}{\pi}
\arccos(\alpha) \approx 0.308799876.$$

\begin{figure}[t]
  \includegraphics[width=1.0\textwidth]{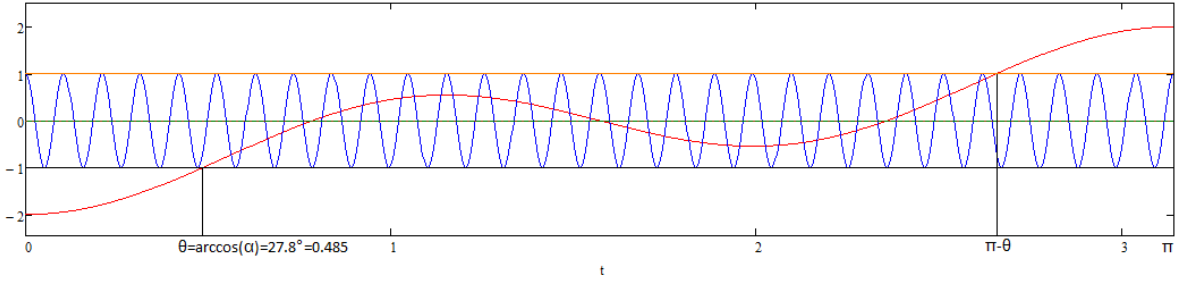}
\caption{
Graph of $f_1=\cos 60 t$ and $f_2=-\cos t-\cos 2t$. Each intersection point of these graphs corresponds to an unimodular root of the reciprocal polynomial $x^{120}+x^{63}+x^{61}+x^{59}+x^{57}+1$.}
\label{fig:2}       
\end{figure}

Each intersection point of graph of $f_1=\cos 60 t$ and $f_2=-\cos t-\cos 2t$ (see fig. 2) corresponds to an unimodular root of the reciprocal polynomial $x^{120}+x^{63}+x^{61}+x^{59}+x^{57}+1$. Nonunimodular roots have arguments in $[-\theta,\theta]$ or in $[\pi-\theta,\pi+\theta]$ where $\theta=\arccos(\alpha)\approx 0.485 \approx 27.8^{\circ}$ (see fig. 3). Since there are 41 intersection point on $[0,\pi]$ and $f_1$, $f_2$ are both even it follows that there are $120-2\cdot41=38$ nonunimodular roots so that $C(P_{120})=38/120\approx 0.317$ is close to the limit of
$C(P_{2n})\rightarrow \frac{2}{\pi}\arccos(\alpha)=\frac{2}{\pi}\theta \approx 0.308799876$. 
\begin{figure}[t]
  \centering
  \includegraphics[width=0.6\textwidth]{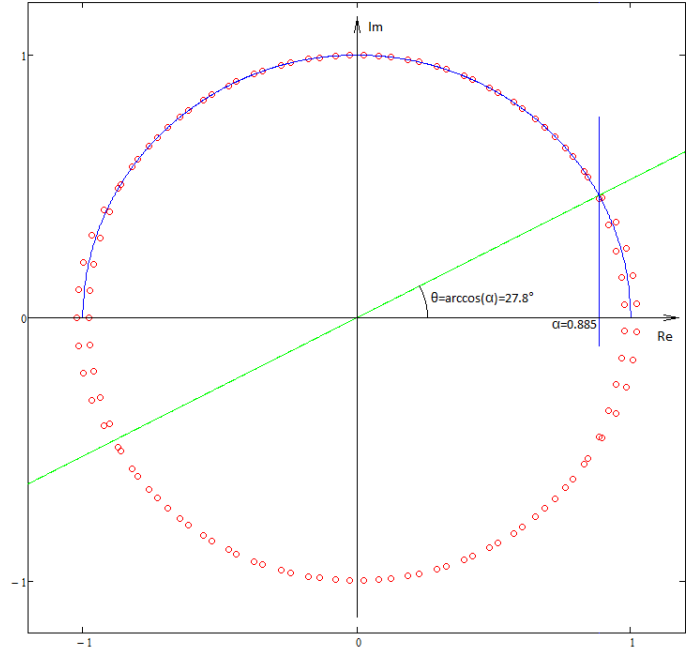}
\caption{Roots of the reciprocal polynomial $x^{120}+x^{63}+x^{61}+x^{59}+x^{57}+1$ are represented with $\circ$. Nonunimodular roots have arguments in $[-\theta,\theta]$ or in $[\pi-\theta,\pi+\theta]$ where $\theta=\arccos(\alpha)\approx 0.485 \approx 27.8^{\circ}$.
}
\label{fig:3}       
\end{figure}

\subsection{Heptanomials with smallest limit points of $C(P_{2n})$}

Between all heptanomials $x^{2n}+a_{j_2}x^{n+j_2}+a_{j_1}x^{n+j_1}+a_{0}x^{n}+a_{j_1}x^{n-j_1} +a_{j_2}x^{n-j_2}+1$ an exhaustive search such that $j_1=1,2,\ldots,10$,  $j_2=j_1+1,j_1+2,\ldots,10$, $a_{0}=\pm1,\pm 2,\ldots \pm 10$, $a_{j_1}=\pm1,\pm 2,\ldots \pm 10$,$a_{j_2}=\pm1,\pm 2,\ldots \pm 10$ suggests that $C(x^{2n}+x^{n+4}+x^{n+2}+x^n +x^{n-2}+x^{n-4}+1)$ has the minimal limit point.

Using the algorithm we solve the equation
$\cos4t+\cos 2t+1/2=\pm1$. If we develop $\cos 4t$ and $\cos2 t$ and substitute $\cos t=x$ we get biquadratic equation $8x^4-6x^2+1/2=1$ with four real solutions $\pm \sqrt{\frac{3}{8}\pm\frac{\sqrt{13}}{8}}$ and $8x^4-6x^2+1/2=-1$ without any real solution. Using symmetry we can show that
$$C(P_{2n})\rightarrow \frac{2}{\pi}
\arccos\left(\sqrt{\frac{3}{8}+\frac{\sqrt{13}}{8}}\right) \approx 0.2741871146.$$

\subsection{Octanomials with smallest limit points of $C(P_{2n})$}

Between all octanomials
$x^{2n}+a_{j_3}x^{n+j_3}+a_{j_2}x^{n+j_2}+a_{j_1}x^{n+j_1}+a_{j_1}x^{n-j_1} +a_{j_2}x^{n-j_2}+a_{j_3}x^{n-j_3}+1$
an exhaustive search such that $j_1=1,2,\ldots,10$,  $j_2=j_1+1,j_1+2,\ldots,10$, $j_3=j_2+1,j_2+2,\ldots,10$,
$a_{j_1}=\pm1,\pm 2,\ldots \pm 10$,$a_{j_3}=\pm1,\pm 2,\ldots \pm 10$, $a_{j_3}=\pm1,\pm 2,\ldots \pm 10$, suggests that $C(x^{2n}+x^{n+5}+x^{n+3}+x^{n+1}+x^{n-1} +x^{n-3}+x^{n-5}+1)$ has the minimal limit point.

Using the algorithm we solve the equation
$\cos5t+\cos 3t+\cos t=\pm1$. If we develop $\cos 5t$ and $\cos 3t$ and substitute $\cos t=x$ we get pentic equation $16x^5-16x^3+3x=\pm 1$ with two real solutions $\pm 0.92757157104393247625$. Using symmetry we can show that
$$C(P_{2n})\rightarrow \frac{2}{\pi}\arccos(0.92757157104393247625) \approx 0.24378469902904315.$$

\subsection{Nonanomials with smallest limit points of $C(P_{2n})$}

Between all nonanomials $x^{2n}+a_{j_3}x^{n+j_3}+a_{j_2}x^{n+j_2}+a_{j_1}x^{n+j_1}+a_{0}x^{n}+a_{j_1}x^{n-j_1} +a_{j_2}x^{n-j_2}+a_{j_3}x^{n-j_3}+1$ an exhaustive search such that $j_1=1,2,\ldots,10$,  $j_{i}=j_{i-1}+1,j_{i-1}+2,\ldots,10$, $i=2,3$; $a_{j}=\pm1,\pm 2,\ldots \pm 10$, $j=0,j_1,j_2,j_3$ suggests that $C(x^{2n}+x^{n+6}+x^{n+4}+x^{n+2}+x^n +x^{n-2}+x^{n-4}+x^{n-6}+1)$ has the minimal limit point.

Using the algorithm we solve the equation
$\cos6t+\cos4t+\cos 2t+1/2=\pm1$. If we develop $\cos6t$, $\cos 4t$ and $\cos2 t$ and substitute $\cos t=x$ we get bicubic equation $32x^6-40x^4+12x^2-1/2=1$ with two real solutions
$$\alpha_{1,2}=\pm\frac{\sqrt{60\sqrt[3]{\frac{\sqrt{29}}{384}+\frac{61}{3456}}+144\sqrt[3]{\frac{61\sqrt{29}}{663552}+\frac{3035}{5971968}}+7}}{12\sqrt[6]{\frac{\sqrt{29}}{384}+\frac{61}{3456}}}$$
and $32x^6-40x^4+12x^2-1/2=-1$ without any real solution. Using symmetry we can show that
$$C(P_{2n})\rightarrow \frac{2}{\pi}
\arccos(|\alpha_1|) \approx 0.21854988117598984.$$

\subsection{Decanomials with smallest limit points of $C(P_{2n})$}

Between all decanomials $x^{2n}+a_{j_4}x^{n+j_4}+a_{j_3}x^{n+j_3}+a_{j_2}x^{n+j_2}+a_{j_1}x^{n+j_1}+a_{j_1}x^{n-j_1} +a_{j_2}x^{n-j_2}+a_{j_3}x^{n-j_3}+a_{j_4}x^{n-j_4}+1$ an exhaustive search such that $j_1=1,2,\ldots,10$,  $j_{i}=j_{i-1}+1,j_{i-1}+2,\ldots,10$, $i=2,3,4$; $a_{j}=\pm1,\pm 2,\ldots \pm 10$, $j=j_1,j_2,j_3,j_4$ suggests that $C(x^{2n}+x^{n+7}+x^{n+5}+x^{n+3}+x^{n+1}+x^{n-1} +x^{n-3}+x^{n-5}+x^{n-7}+1)$ has the minimal limit point.

Using the algorithm we solve the equation
$\cos7t+\cos5t+\cos 3t+\cos t=\pm1$. If we develop $\cos7t$, $\cos 5t$ and $\cos3 t$ and substitute $\cos t=x$ we get two equations $64x^7-96x^5+40x^3-4x=\pm1$ with two real solutions
$\pm 0.9521755884525 $. Using symmetry we can show that
$$C(P_{2n})\rightarrow \frac{2}{\pi}\arccos(0.95217558845251615756) \approx 0.19768155115418617.$$

We remark that this is the smallest limit point of $C(P_{2n})$, $n\rightarrow\infty$ we know.

\subsection{Polynomials with smallest limit points of $C(P_{2n})$}

Our calculations suggest that the next conjecture seems to be true:
\begin{conjecture}
If $P_{2n}$ is a sum of $2k+3$ monomials, i.e. $a_0\ne 0$, then the sequence $C(x^{2n}+x^{n+2k}+\cdots+x^{n+4}+x^{n+2}+x^n +x^{n-2}+x^{n-4}+\cdots+x^{n-2k}+1)$ tends to the smallest limit, greater than zero, of $C(P_{2n})$, $n\rightarrow\infty$.

Let $k\ge 1$ be an integer. If $P_{2n}$ is a sum of $2k+4$ monomials, i.e. $a_0=0$, then the sequence $C(x^{2n}+x^{n+2k+1}+\cdots+x^{n+5}+x^{n+3}+x^{n+1}+x^{n-1} +x^{n-3}+x^{n-5}+\cdots+x^{n-2k-1}+1)$ tends to the smallest limit, greater than zero, of $C(P_{2n})$, $n\rightarrow\infty$.
\end{conjecture}

But in the case of dodecanomials we found that $C(x^{2n}+x^{n+9}+x^{n+7}+2x^{n+5}+2x^{n+3}+2x^{n+1}+2x^{n-1} +2x^{n-3}+2x^{n-5}+x^{n-7}+x^{n-9}+1)$ tends to $2\arccos(0.943468)/\pi=0.215085$ which is smaller than $0.226163=2(\arccos(0.966357)+\arccos(0.877575)-\arccos(0.919147))/\pi$ the limit of $C(x^{2n}+x^{n+9}+x^{n+7}+x^{n+5}+x^{n+3}+x^{n+1}+x^{n-1} +x^{n-3}+x^{n-5}+x^{n-7}+x^{n-9}+1)$.
Nevertheless the conjecture seems to be true for many $k$.

$C(x^{2n}+x^{n+9}+x^{n+7}+2x^{n+5}+2x^{n+3}+2x^{n+1}+2x^{n-1} +2x^{n-3}+2x^{n-5}+x^{n-7}+x^{n-9}+1)$ $\rightarrow$ $2\arccos(0.943468)/\pi$ $=0.215085<$ $0.226163=$ $2(\arccos(0.966357)+\arccos(0.877575)-\arccos(0.919147))/\pi$ the limit of $C(x^{2n}+x^{n+9}+x^{n+7}+x^{n+5}+x^{n+3}+x^{n+1}+x^{n-1} +x^{n-3}+x^{n-5}+x^{n-7}+x^{n-9}+1)$.

It is natural to ask: do exist
$$\lim_{k\rightarrow\infty}\lim_{n\rightarrow\infty}C\left(x^n\left(x^n+\frac{1}{x^n}+\sum
_{j=1}^k \left(x^{2j-1}+\frac{1}{x^{2j-1}}\right)\right)\right),$$
$$\lim_{k\rightarrow\infty}\lim_{n\rightarrow\infty}C\left(x^n\left(x^n+1+\frac{1}{x^n}+\sum
_{j=1}^k \left(x^{2j}+\frac{1}{x^{2j}}\right)\right)\right)?$$

It is an easy exercise to prove that
$$\cos t+\cos3t+\cdots+\cos(2m-1)t=\frac{\sin2mt}{2\sin t}  ,$$
$$\frac{1}{2}+\cos 2t+\cos4t+\cdots+\cos2mt=\frac{\sin(2m+1)t}{2\sin t}  .$$

These formulae enable us to calculate $f_2$ and $C(P_{2n})$ much faster. Our experiments with $k\approx $ half of million, $n\approx $ one hundred million suggest that these limits exist and that they are both equal to $0.20885$.

\section{Extension of Mahler measure}

The definition of the Mahler measure could be extended to polynomials in several
variables. We recall Jensen's formula which states that
$\frac{1}{2\pi}\int_0^{2\pi} \log |P(e^{ i\theta})|\d \theta= \log |a_0| + \sum_{j=1}^d \log\max (|\alpha_j|, 1)$
Thus
$$M(P) = \exp \left\{ \frac{1}{2\pi}\int_0^{2\pi} \log |P(e^{ i\theta})|\d \theta\right\},$$
so M(P) is just the geometric mean of $|P(z)|$ on the torus $T$.

Hence a natural
candidate for $M(F)$ is
$$M(F) = \exp\left\{\frac{1}{(2\pi)^r}\int_0^{2\pi} \d \theta_1\cdots\int_0^{2\pi}\log |F(e^{ i\theta_1},\ldots,e^{ i\theta_r})|\d \theta_r \right\}.$$

The smallest known Mahler measures in two variables are
$$M((x + 1)y^2 + (x^2 + x + 1)y + x(x + 1)) = 1.25542\ldots$$ and
$$M(y^2 + (x^2 + x + 1)y + x^2) = 1.28573\ldots$$

Boyd proved (1981)\cite{Boy1} the next
\begin{theorem}\label{ThBoyd}
As $m\rightarrow \infty$, $M(P (x, x^m))\rightarrow M(P(x,y))$.
\end{theorem}

Let $Q(x_1,x_2,\ldots,x_r)=\sum_{j=1}^{k}a_jx_1^{e_{j1}}x_2^{e_{j2}}\cdots x_r^{e_{jr}}$, $a_j\in \mathbb{R}$,
$e_{ji}\in \mathbb{Z}$, $j=1,2,\ldots,k$, $i=1,2,\ldots,r$ and let $W(x_1,x_2,\ldots,x_r)=$
$$=(x_1^nx_2^n\cdots x_r^n)\left(x_1^nx_2^n\cdots x_r^n+x_1^{-n}x_2^{-n}\cdots x_r^{-n}+Q(x_1,x_2,\ldots,x_r)+Q(x_1^{-1},x_2^{-1},\ldots,x_r^{-1})\right).$$
where $n\in \mathbb{N}$ is greater than $\max(|e_{ji}|)$ so that $W$ is a multivariate polynomial. Let $g_2(x_1,x_2,\ldots,x_r):=-1/2(Q(x_1,x_2,\ldots,x_r)+Q(x_1^{-1},x_2^{-1},\ldots,x_r^{-1}))$ and
\begin{equation}\label{eq:g2}
    \overline{g_2(x_1,x_2,\ldots,x_r)} := \begin{cases}
               1, & |g_2(x_1,x_2,\ldots,x_r)|\ge 1\\
               0, & \text{otherwise},
           \end{cases}
\end{equation}
then we can define
\begin{equation}\label{eq:LCint}
LC(W):=\frac{1}{(2\pi)^r}\int_{0}^{2\pi}\int_{0}^{2\pi}\cdots \int_{0}^{2\pi}\overline{g_2(\exp(it_1),\exp(it_2),\ldots,\exp(it_r))}\d t_1\d t_2\cdots \d t_r.
\end{equation}
If $r=1$, $Q(x_1)=a_0/2+\sum_{j=1}^k a_jx_1^j$ then $W(x_1)=P_{2n}(x_1)$, $g_2(\exp(it))=f_2(t)$ and $\overline{g_2(\exp(it))}=\overline{f_2(t)}$ so that, recalling \eqref{eq:Cint},
we conclude that $LC(P_{2n})=\lim_{n\rightarrow\infty}C(P_{2n})$.

If $Q(x,y)=x+y+1$ we can prove that the Boyd's property for $LC$ is valid: $LC(W (x, x^m))\rightarrow LC(W(x,y))$ as $m\rightarrow \infty$. Indeed
$LC(x^{2n}+x^n(x+ x^m+2+x^{-1}+x^{-m}))+1)=\lim_{n\rightarrow\infty}C(x^{2n}+x^n(x+ x^m+2+x^{-1}+x^{-m}))+1)=\frac{1}{2\pi}\int_{0}^{2\pi}\overline{f_{2m}(t)}$, where $f_{2m}(t)=-1-\cos(mt)-\cos(t)$.
Since $\cos(t)=-\cos(\pi-t)$ it follows that $\cos(2(m_1+1)t)=-\cos(\pi-2(m_1+1)t)=-\cos((2m_1+1)(\pi-t))$. Therefore if $m$ is odd then for each interval $I=[a,b]\subseteq [0,\pi]$ such that $|f_{2m}(t)|>1$, $a<t<b$, there is the interval $I'=[\pi-b,\pi-a]$ of the equal length such that $|f_{2m}(t)|\le 1$, $t\in I'$. We conclude that $\frac{1}{\pi}\int_{0}^{\pi}\overline{f_{2m}(t)}=0.5$ for $m$ odd so that $LC(W (x, x^m))\rightarrow 0.5$ as $m\rightarrow \infty$.

On the other hand
\begin{align*}
  LC(W) & =\frac{1}{(2\pi)^2}\int_{0}^{2\pi}\int_{0}^{2\pi}\overline{-\frac{1}{2}(\exp(it_1)+\exp(it_2)+2+\exp(-it_1)+\exp(-it_2))}\d t_1\d t_2 \\
   & =\frac{1}{(2\pi)^2}\int_{0}^{2\pi}\int_{0}^{2\pi}\overline{(-1-\cos t_1-\cos t_2)}\d t_1\d t_2.
\end{align*}
Since if $t_1,t_2\in [0,\pi]$ then $|1+\cos t_1+\cos t_2|\ge 1$ is equivalent with $t_2\le \pi-t_1$ and using the symmetry of the set $\{(t_1,t_2)\in [0,2\pi]\times [0,2\pi]:\overline{-1-\cos t_1-\cos t_2}\ge 1\}$ it follows that
$$LC(W)=\frac{4}{(2\pi)^2}\int_{0}^{\pi}\left(\int_{0}^{\pi-t_1}\d t_2\right)\d t_1=\frac{1}{\pi^2}\int_{0}^{\pi}\left(\pi-t_1\right)\d t_1=\frac{1}{2}.$$
This example as well as numerical approximations of $LC$ of many other polynomials in two variables using the formula \eqref{eq:Cint} and the definition \eqref{eq:LCint} suggest us that the Boyd's limit formula in Theorem \ref{ThBoyd} is also valid for $LC$ i.e. we propose the following
\begin{theorem}\label{ThJJ}
As $m\rightarrow \infty$, $LC(W (x, x^m))\rightarrow LC(W(x,y))$.
\end{theorem}

To prove Theorem \ref{ThJJ} we use two lemmas which Everest and Ward proved in \cite{EW}.
For the sake of completeness we cite proofs of the lemmas in the Appendix. Denote as usual the (multiplicative) circle group by $K=\mathbb{S}^1$, and the torus by $K^2=\mathbb{S}^1\times\mathbb{S}^1$. For an integrable function $f:K\rightarrow \mathbb{C}$, write
\begin{equation}\label{eq:muK}
\int_{0}^{1}f(e^{2\pi i \theta})d\theta=\int f(x)d\mu_K= \int fd\mu_K
\end{equation}

and for an integrable function $g:K^2\rightarrow \mathbb{C}$, write
\begin{equation}\label{eq:muK2}
\int_{0}^{1}\int_{0}^{1}g(e^{2\pi i \theta_1},e^{2\pi i \theta_2})d\theta_1d\theta_2=\int g(x_1,x_2)d\mu_{K^2}= \int gd\mu_{K^2}.
\end{equation}
We will use the Lebesgue measure $\mu_K$ on the circle to evaluate the measure of disjoint unions of intervals (whose measure is simply the sum of the lengths).

\begin{lemma}\label{L3t22}
Let $\phi:K^2\rightarrow \mathbb{R}$ be any continuous function. Then
$$\lim_{N\rightarrow\infty}\int \phi(x,x^N)d\mu_{K}=\int \phi d\mu_{K^2}$$
\end{lemma}

\begin{lemma}\label{LExer}
Let $\phi:K^2\rightarrow \mathbb{R}$ be any Riemann-integrable function and $\delta>0$ be given. There are finite trigonometric series $P(x_1,x_2)=\sum_{||\textbf{n}||<M}a_{\textbf{n}}x_1^{n_1}x_2^{n_2}$ and $Q(x_1,x_2)=\sum_{||\textbf{n}||<M}b_{\textbf{n}}x_1^{n_1}x_2^{n_2}$ with the property that
$$P(x_1,x_2)\le \phi(x_1,x_2)\le Q(x_1,x_2)$$
for all $(x_1,x_2)\in K^2$ and
$$\int(Q-P)d\mu_{K^2}<\delta,$$ 
where $||\textbf{n}||=\max\{|n_1|,|n_2|\}$, $M>0$.
\end{lemma}
The finite sums $P$, $Q$ used to bound $\phi$ are called trigonometric polynomials since in the additive group notation the monomial $x_1^{n_1}x_2^{n_2}$  corresponds to $e^{2\pi i(n_1\theta_1+n_2\theta_2)}$ under the correspondence $x_1=e^{2\pi i \theta_1}$, $x_2=e^{2\pi i \theta_2}$.

\noindent\textit{Proof of Theorem \ref{ThJJ}}.
If $r=2$ in \eqref{eq:g2} then function $\overline{g_2(x_1,x_2)}$
is not continuous but is Riemann-integrable (since it is bounded and the set of discontinuities
of $g_2$ has measure 0). By Lemma \ref{LExer} there are finite trigonometric series $P_2(x_1,x_2)=\sum_{||\textbf{n}||<M}a_{\textbf{n}}x_1^{n_1}x_2^{n_2}$ and $Q_2(x_1,x_2)=\sum_{||\textbf{n}||<M}b_{\textbf{n}}x_1^{n_1}x_2^{n_2}$ with the property that
$$P_2(x_1,x_2)\le \overline{g_2(x_1,x_2)}\le Q_2(x_1,x_2)$$
for all $(x_1,x_2)\in K^2$ and
$$\int(Q_2-P_2)d\mu_{K^2}<\delta.$$
It follows that
\begin{equation}\label{eq:PgQ}
\int P_2(x,x^m)d\mu_K\le \int\overline{g_2(x,x^m)}d\mu_K\le \int Q_2(x,x^m)d\mu_K.
\end{equation}
Function $P_2$, $Q_2$ are continuous so by Lemma \ref{L3t22}
\begin{equation}\label{eq:PtP}
\int P_2(x,x^m)d\mu_K\rightarrow \int P_2(x_1,x_2)d\mu_{K^2},
\end{equation}
\begin{equation}\label{eq:QtQ}
\int Q_2(x,x^m)d\mu_K\rightarrow \int Q_2(x_1,x_2)d\mu_{K^2}.
\end{equation}
Since $\int Q_2d\mu_{K^2}-\int P_2d\mu_{K^2}<\delta$ and $\delta>0$ was arbitrary, \eqref{eq:PgQ}, \eqref{eq:PtP} and \eqref{eq:QtQ} then show that

$$\int\overline{g_2(x,x^m)}d\mu_K\rightarrow \int \overline{g_2(x_1,x_2)}d\mu_{K^2}.$$
Recalling \eqref{eq:muK}, \eqref{eq:muK2} and using the substitutions $t=2\pi\theta$, $t_1=2\pi\theta_1$, $t_2=2\pi\theta_2$ it follows that
$$\int g_2(x,x^m)d\mu_K=\int_0^1 g_2(e^{2\pi i\theta},e^{2m\pi i\theta})d\theta=\frac{1}{2\pi}\int_0^{2\pi} g_2(e^{it},e^{m it})dt=LC(W(x,x^m))$$
and
$$\int \overline{g_2(x_1,x_2)}d\mu_{K^2}=\int_0^1\int_0^1 g_2(e^{2\pi i\theta_1},e^{2\pi i\theta_2})d\theta_1 d\theta_2=\frac{1}{4\pi^2}\int_0^{2\pi}\int_0^{2\pi} g_2(e^{it_1},e^{it_2})dt_1dt_2=LC(W(x_1,x_2)),$$
hence the claim follows.


\textbf{Acknowledgements:} I am grateful to Jonas Jankauskas for his careful reading and for sending me a sketch of a proof of Theorem \ref{ThJJ}.

\section{Appendix}
\subsection{Proof of Lemma \ref{L3t22}}
By the Stone-Weierstrass approximation theorem (see \cite{EW} Appendix B) for any $\epsilon>0$ there is an $M$ and there are coefficients $\{a_{\textbf{n}}\}_{||\textbf{n}||<M}$ for which
$$\left|\phi(x_1,x_2)-\sum_{||\textbf{n}||<M}a_{\textbf{n}} x_1^{n_1}x_2^{n_2}\right| <\epsilon$$
for all $x_1,x_2 \in K$, where $||\textbf{n}||=\max\{|n_1|,|n_2|\}$. Since this estimate is uniform, it is enough to show that the required convergence happens for the function
$$G(x_1^{n_1}x_2^{n_2})=\sum_{||\textbf{n}||<M}a_{\textbf{n}} x_1^{n_1}x_2^{n_2}.$$
Now
$$\int G d \mu_{K^2}= \sum_{||\textbf{n}||<M} \int a_{\textbf{n}} x_1^{n_1}x_2^{n_2} d \mu_{K^2}=a_0.$$
On the other hand
$$\int G(x,x^N) d \mu_{K}= \sum_{||\textbf{n}||<M} \int a_{\textbf{n}} x_1^{n_1+Nn_2} d \mu_{K}=\sum_{\textbf{n}:n_1+Nn_2=0} a_{\textbf{n}}.$$
For fixed $M$, when $m$ is large $n_1+Nn_2=0$ with $||\textbf{n}||<M$ if and only if $n_1=n_2=0$, so
$$\lim_{m\rightarrow\infty}\int G(x,x^N) d \mu_{K}=a_0=\int G d \mu_{K^2}  .$$

\subsection{Proof of Lemma \ref{LExer}}

By the Stone-Weierstrass approximation theorem (see \cite{EW} Appendix B), it is enough to find continuous functions $P$ and $Q$ with the stated properties. Since $\phi$ is Riemann -integrable it is bounded above by $R$ say, and there is a finite collection of rectangles $\{A_1,A_2,\ldots,A_n\}$ with the property that
$$\sum_{i=1}^{n} \textrm{area} A_i \cdot \sup_{(x_1,x_2)\in A_i}\{\phi(x_1,x_2)\}-\int \phi d \mu_{K^2}<\delta/4.$$
The word "rectangle" means a set of the form
$$\{ (e^{2\pi i \theta_1},e^{2\pi i \theta_2}):\theta_j\in[a_j,b_j]]\subset \mathbb{T}.\}$$
Now define

\begin{equation*}
    \overline{Q}(x_1,x_2) = \begin{cases}
               \sup_{(x_1,x_2)\in A_i}\{\phi(x_1,x_2)\} & \textrm{if} (x_1,x_2)\in A_i\setminus\bigcup_{j\ne i} A_j,\\
               R & \text{otherwise.}
           \end{cases}
\end{equation*}

Then $ \overline{Q}\ge \phi$ and their integrals are within $\delta/4$. Now approximate $ \overline{Q}$ from above by a continuous function $Q$ which is equal to $ \overline{Q}$ except very close to the boundary of each rectangle $A_i$, and fills in continuously to reach the value $R$ on the boundary. This can be done by keeping $\int Qd\mu_{K^2} -\int \overline{Q}d\mu_{K^2}<\delta/4$, which gives a continuous $Q$ with $Q\ge \phi$ and $\int Qd\mu_{K^2} -\int \phi d\mu_{K^2}<\delta/2$.

Repeating the argument from below (or simply repeating the argument for $-\phi$) gives $P$.

\end{document}